\documentclass[11pt, lettersize, notitlepage, leqno, footskip=0.6cm]{article} 
\usepackage{amssymb}
\usepackage{amsthm}
\usepackage{amsmath}	
\usepackage{graphicx}
\usepackage{amscd}

\usepackage{thmtools}
\usepackage[none]{hyphenat}	

\newcommand{\bc}{\mathbb C}
\newcommand{\bF}{\mathbb F}

\newcommand{\bz}{\mathbb Z}

\newcommand{\bq}{\mathbb Q}
\newcommand{\br}{\mathbb R}

\newcommand{\bT}{\mathbf T}

\newcommand{\Gal}{\mathrm{Gal}}

\newcommand{\W}{\mathcal W}
\newcommand{\WA}{\mathcal{W}_{A_v}}

\newcommand{\zbar}{\overline {\mathbb{Z}}}
\newcommand{\qbar}{\overline {\mathbb{Q}}}

\newcommand{\Ebar}{\overline {E}}

\newcommand{\bg}{\mathbb{G}}
\newcommand{\ev}{\mathrm{ev}}

\newcommand{\absq}{\mathrm{Gal}_{\bq}}
\newcommand{\absql}{\mathrm{Gal}_{\bq_l}}
\newcommand{\absqp}{\mathrm{Gal}_{\bq_p}}

\newcommand{\absf}{\mathrm{Gal}_F}

\newcommand{\abse}{\mathrm{Gal}_E}
\newcommand{\absk}{\mathrm{Gal}_K}

\newcommand{\Fr}{\mathrm{Fr}}

\newcommand{\la}{\langle}
\newcommand{\ra}{\rangle}

\newcommand{\lra}{\longrightarrow}
\newcommand{\xra}{\xrightarrow}
\newcommand{\hra}{\hookrightarrow}

\newcommand{\bs}{\backslash}

\newcommand{\wti}{\widetilde}
\newcommand{\mf}{\mathfrak}
\newcommand{\mc}{\mathcal}
\newcommand{\mb}{\mathbb}
 
\newcommand{\mr}{\mathrm}

\newcommand{\mfp}{\mathfrak{p}}

\newcommand{\mfl}{\mathfrak{l}}

\newcommand{\Br}{\mr{Br}}

\newcommand{\al}{\alpha}
\newcommand{\be}{\beta}

\newcommand{\Gam}{\Gamma}

\newcommand{\lamb}{\lambda}

\newcommand{\si}{\sigma}

\newcommand{\om}{\omega}
\newcommand{\Om}{\Omega}

\newcommand{\ov}{\overline}

\newcommand{\sub}{\subseteq}
\newcommand{\nsub}{\nsubseteq}

\newcommand{\bqW}{\bq(\W_{A_v})}

\newcommand{\GalW}{\Gal(\bq(\W_{A_v})/\bq)}

\newcommand{\ovm}{\ov{\om}}

\DeclareMathOperator{\MT}{MT}

\DeclareMathOperator{\Aut}{Aut}
\DeclareMathOperator{\Cent}{Cent}

\DeclareMathOperator{\End}{End}
\DeclareMathOperator{\GL}{GL}

\DeclareMathOperator{\Hom}{Hom}
\DeclareMathOperator{\Cor}{Cor}

\DeclareMathOperator{\Mat}{Mat}

\DeclareMathOperator{\SL}{SL}

\DeclareMathOperator{\rank}{rank\;}

\DeclareMathOperator{\inv}{inv}

\DeclareMathOperator{\der}{der}
\DeclareMathOperator{\Nm}{Nm}
\DeclareMathOperator{\charac}{char}
\DeclareMathOperator{\Res}{Res}

\DeclareMathOperator{\conn}{conn}

\DeclareMathOperator{\ord}{ord}

\DeclareMathOperator{\IV}{IV}

\newcommand{\Conj}{\mathbf{Conj}}

\newcommand{\SA}{\mc{S}_A}

\newcommand{\muN}{\eta_{_{0,v}}}
\newcommand{\muH}{\mu_{_H}}

\newcommand{\PAv}{P_{A_{v}}(X)}

\newcommand{\mcO}{\mathcal{O}}

\newtheorem{Thm}{Theorem}[section]
\newtheorem{Prop}[Thm]{Proposition}
\newtheorem{Lem}[Thm]{Lemma}
\newtheorem{Corr}[Thm]{Corollary}

\newtheorem{Def}{Definition}[section]

\declaretheoremstyle[%
  spaceabove=-2pt,%
  spacebelow=8pt,%
  headfont=\normalfont\itshape,%
  postheadspace=1em,%
  qed=\qedsymbol%
]{mystyle} 
\declaretheorem[name={Proof},style=mystyle,unnumbered,
]{prf}

\numberwithin{equation}{section}

\title{On some abelian varieties of type $\IV$}
\author{Steve Thakur}  
\date{\vspace{-3ex}}

\begin{document} 

\maketitle

\begin{abstract} We study a certain class of simple abelian varieties of type $\IV$ (in Albert's classification) over number fields with Mumford-Tate groups of type $A$. In particular, we show that such abelian varieties have ordinary reduction away from a set of places of Dirichlet density zero, thus confirming a special case of a broader conjecture of Serre's. We also study the splitting types and Newton polygons of the reductions of the abelian varieties of this type with small dimension (nine).
\end{abstract}

\begin{center}
\section{\fontsize{11}{11}\selectfont Introduction}
\end{center}
\vspace{-0.1cm}

For an absolutely simple abelian variety $A$ over a number field $F$, it is natural to study the reduction $A_v$ at each place $v$ of good reduction. In particular, two questions of significant interest are:\\ 
1. to determine the splitting of $A_v$ into simple abelian varieties up to isogeny\\ 
2. to determine the Newton polygon of $A_v$.

The following well known conjectures suggest that, after passing to a suitable finite extension, such an abelian variety has iso-simple and ordinary reduction at a set of places of density one.

\noindent \textbf{Conjecture 1.} ([MP08], [Zyw14]) \textit{Let $X$ be an absolutely simple abelian variety over a number field $F$. After replacing $F$ by a finite extension if necessary, there exists a density one set of places such that the reduction $X_v$ is isogenous to the $d$-th power of a simple abelian variety where $d^2$ is the dimension of the division algebra $\End^0_{\qbar}(X)$ over its center.}

\noindent \textbf{Conjecture 2.} ([Pin98]) \textit{Let $X$ be an absolutely simple abelian variety over a number field $F$. After replacing $F$ by a finite extension if necessary, there exists a density one set of places such that the reduction $X_v$ is ordinary.}
 
In this article, we study a certain class of abelian varieties that are of type $\IV$ in Albert's classification and have Mumford-Tate groups (and $l$-adic monodromy groups) of small rank and of type $A$. In particular, we verify Conjecture 2 for this class of abelian varieties, relying heavily on the techniques of [Pin98]. The following theorem is the main result of this article.

\begin{Thm} Let $A$ be an abelian variety of type $\IV(1,d,\wedge^r(\mr{Std}))$ over a number field. If the $l$-adic monodromy groups $G_{A,l}$ are connected, then $A$ has ordinary reduction at a set of places of Dirichlet density one.\end{Thm} 

\noindent The notation in the statement of the theorem is explained in Definition 1.1.

\vspace{-0.1cm}
\subsection{\fontsize{11}{11}\selectfont  Notations and Terminology}
\vspace{-0.2cm }
For an abelian variety $A$ over a number field and a prime $l$, the \textit{$l$-adic monodromy subgroup} $G_{A,l}$ is the Zariski closure of the image of $\rho_{_l}:\absf\lra \GL(V_l(A))$ and $G_l^0(A)$ is the connected component of the identity. The image $\rho_l(\absf)$ is a compact $l$-adic analytic subgroup of $G_l(\bq_l)$ which is of finite index in $G_l(\bq_l)$ ([Bog80]). We denote the Zariski closure of $\rho_{_l}(\absf)$ in $\GL(T_l(A))$ by $\mc{G}_l$. It is a group scheme over $\bz_l$ with generic fiber $G_l$.

For a finite place $v$ of $F$ with $\charac(v) =p$, the Zariski closure of the image $\rho_{_p}(D_v)$ of a decomposition group $D_v\sub \absf$ under the map $\rho_{_p}:\absf\lra \GL(V_p(A))$ is denoted by $H_v(A)$. 

Let $g$ be the dimension of $A$. The Betti cohomology group $V(A)=H^1(A(\bc),\bq)$ is a $2g$-dimensional vector space over $\bq$ endowed with a decomposition $V(A)\otimes_{\bq}\bc = V^{1,0}\oplus V^{0,1}$ such that $\overline{V^{1,0}} = V^{0,1}$. Let $\mu_{\infty}:\bg_{m,\bc}\lra \GL_{2g,\bc}$ be the cocharacter through which any $z\in \bc^{*}$ acts by multiplication by $z$ on $V^{1,0}$ and trivially on $V^{0,1}$. The \textit{Mumford-Tate group} $\MT_A$ is the unique smallest $\bq$-algebraic subgroup of $\GL_{2g,\bq}$ such that $\mu_{\infty}$ factors through $\MT_A\times_{\bq} \bc$.

The following long-standing conjecture suggests an intrinsic relation between the two notions:\\

\noindent\textbf{Conjecture.} (Mumford-Tate) \textit{For an abelian variety $A$ over a number field and a prime $l$, $\MT_A\times_{\bq}\bq_l = G_l^0(A)$}.

\vspace{-0.2cm}
\subsection{\fontsize{11}{11}\selectfont  Some background}
\vspace{-0.1cm}

For an abelian variety $A$ over field $F$ of characteristic zero, we denote by $\End^0(A)$ the endomorphism algebra $\End_{\qbar}(A)\otimes_{\bz}\bq$. According to Albert's classification, a simple abelian variety $A$ of dimension $g$ has one of the following possibilities for its endomorphism algebra:

\noindent \underline{Type $\mr{I}$}: $\End^0(A)$ is a totally real field of degree $e$ such that $e|g$.\\
\noindent \underline{Type $\mr{II}$}: $\End^0(A)$ is a totally indefinite quaternion algebra over a totally real field of degree $e$ with $2e|g$.\\
\noindent \underline{Type $\mr{III}$}: $\End^0(A)$ is a totally definite quaternion algebra over a totally real field of degree $e$ with $2e|g$.\\
\noindent \underline{Type $\mr{IV}$}: $D:=\End^0(A)$ is a central division algebra of dimension $d^2$ over a CM field $K$ of degree $2e$ with $ed^2|2g$ and is equipped with an involution $*$ of the second kind.

In this article, we will be studying a certain class of abelian varieties of type $\IV$ whose Mumford-Tate groups are of type $A$. Note that if $A$ is such an abelian variety of dimension $g$ and $D$ its endomorphism algebra, the Hasse invariants of $D$ are such that for any place $v$ of $K$,\\ $\inv_v(D)+\inv_{\ov{v}}(D)=0$.

Abelian varieties of type $\IV$ such that $D=K$ have been studied in [Chi92]. Our focus will be a certain class for which $D$ is non-commutative, which seems to be a relatively less explored class of abelian varieties for questions of this type. 

We now state a few fundamental theorems we shall need repeatedly in this article. For an abelian variety over $A$ over a number field $F$, we define the field $$F_A^{\conn}:=\bigcap\limits_l F(A[l^{\infty}]). $$ This is the minimum extension of $F$ such that the $l$-adic monodromy groups are connected and is a number field Galois over $\bq$ ([LP95]).

\begin{Thm} $\mr{(Faltings)}\;$ Let $X$ be an abelian variety over a number field $F$ such that $F_A^{\conn} = F$. Then we have the following:\\
\noindent $\mr{(i)}$ The centralizer of $G_l$ in $\End_{\bq_l}(V_l)$ is $\End(A)\otimes_{\bz} \bq_l$.\\
\noindent $\mr{(ii)}$ The group $G_l$ is reductive.\\
\noindent $\mr{(iii)}$ $\End^0(A_{F}) = \End^0(A_{\qbar})$\end{Thm}
\vspace{-0.1cm}
So, for an abelian variety $A$ with $\End^0(A)$ of dimension $d^2$ over its center, then the $l$-adic representation $V_l(A)$ is the $d$-th power of an absolutely simple representation.

The following theorems of Deligne and Serre give one of the inclusions for the Mumford-Tate conjecture and reduce it to the equalities of the ranks of $\MT_A$ and $G_{A,l}$.

\begin{Thm} $\mr{(Deligne)}\;$  For any prime $l$, $G_l^0\times_{\bq} \bq_l\subseteq \MT_A\times_{\bq} \bq_l$.\end{Thm}

\vspace{-0.3cm}

\begin{Thm} $\mathrm{(Serre)}\;$ For any abelian variety $A$, the $l$-adic monodromy groups $G_{A,l}$ are all of the same rank.\end{Thm}
\vspace{-0.2cm}
Henceforth, we shall refer to this common rank as \textit{the rank of the abelian variety}. Deligne's theorem on the first inclusion combined with Serre's rank independence theorem show that the Mumford-Tate conjecture reduces to the rank of the abelian variety being equal to the rank of the Mumford-Tate group. More precisely:

\begin{Thm} $\mathrm{(Serre)}\;$ If the inequality \vspace{-0.15cm}$$\rank G_p^{0}\leq  \rank \MT_A$$ is an equality for any prime $p$, the Mumford-Tate conjecture is true for $A$.  \end{Thm}
\vspace{-0.1cm}

\noindent\underline{Strict Compatibility.} The $l$-adic representations attached to an abelian variety are \textit{strictly compatible} in the sense of Serre. A proof may be found in [Del74]. 

\begin{Thm} $\mr{[Del74]}$ For an abelian variety $A$, fix a finite set $\mathcal{S}$ of non-archimedean places such that $A$ has good reduction outside $\mathcal{S}$. Let $v\notin \mathcal{S}$, $l\neq \mathrm{char}(k_v)$\\
$1.$ $\rho_{_l}$ is unramified at $v$\\ 
$2.$ The characteristic polynomial of $\rho_{_l}(\mathrm{Fr}_v)$ has coefficients in $\bz$ and is independent of $l$.\end{Thm} 
\vspace{-0.1cm}
\noindent \textbf{Notations:} An immediate consequence is that the eigenvalues are in $\zbar$ and are independent of $l$. We denote this polynomial by $P_{A_v}(X)\in\bz[X]$, its zeros by $\WA$ and the multiplicative group they generate by $\Phi_{A_v}$. The group $\Phi_{A_v}$ is abelian with $\rank \Phi_{A_v}\leq \rank A$. Furthermore, if the field is enlarged to ensure that the $l$-adic monodromy groups are connected, this inequality is an equality for a set of places of Dirichlet density one ([LP95]). 

\begin{Thm}$\mr{[Pin98]}$ Each simple factor of the root system of $G_l^0$ has type $A$, $B$, $C$ or $D$ and its highest weights in the tautological representation are minuscule.\end{Thm}
\vspace{-0.1cm}

\noindent (A \textit{minuscule} representation is a representation such that the Weyl group acts transitively on the set of weights).

In particular, if a simple factor is of type $A_n$, then the corresponding representation must be the $r$-th exterior power of the standard representation of $\SL_{n+1}$ for some $r\leq n$. We now describe the abelian varieties we intend to study in this article.

\begin{Def} \normalfont For an odd integer $d$ and an integer $r\leq d-1$ prime to $d$, we say an abelian variety $A$ over a field of characteristic zero is of type $\IV(1,d,\wedge^r(\mr{Std}))$ if:\\
	$\bullet\; D=\End^0(A)$ is a $d^2$-dimensional central division algebra over an imaginary quadratic field $K$, endowed with an involution $*$ of the second kind.\\
	$\bullet$ The Mumford-Tate group $\MT_A = G_{\wti{D}}/\mu_r$ where $\wti{D}$ is the unique central division algebra over $K$ such that $[\wti{D}]^r = [D]^{-1}$ in $\Br(K)$ and $G_{\wti{D}}(R) := \{\al\in (\wti{D}\otimes_{\bq}R)^{\times}: \al\al^*\in R\}$ for any $\bq$-algebra $R$.\\
	$\bullet\; H^1(A(\bc),\bc)$ is the sum of $2d$ copies of the $r$-the exterior power of the standard representation of $\SL_{d,\bc}$.\end{Def}

We note that a central division algebra over an imaginary quadratic field $K$ is endowed with an involution of the second kind if and only if the corestriction $\Cor_{K/\bq}(D) = 0$ in $\Br(\bq)$ ([KMRT98]).

Abelian varieties over $\bc$ satisfying the conditions in Definition 1.1 are constructed in [Orr15] for the case where $d=2r+1$. We provide a virtually identical construction in the appendix for the case where $\gcd(r,d)=1$. The existence of such abelian varieties over number fields then follows from the following theorem of Noot.

\begin{Thm}$ \mr{([No95],\;Theorem\; 1.7)}$ Let $X$ be an abelian variety over $\bc$. Then there exists an abelian variety $A$ over $\qbar$ such that $\MT_A = \MT_X$, $\End^0(A) = \End^0(X)$ and $A$ satisfies the Mumford-Tate conjecture.\end{Thm}

Such an abelian variety is absolutely simple of dimension $d{d\choose r}$ and satisfies the Mumford-Tate conjecture. The following argument is sketched in [Orr15]. We provide it for the reader's convenience.

\begin{Prop} Suppose either $d>3$ or $r>1$. Then an abelian variety $A$ of type $\IV(1,d,\wedge^r(\mr{Std}))$ satisfies the Mumford-Tate conjecture.\end{Prop}

\begin{prf} Note that for any rational prime $l$, the pair $(\MT_A\otimes_{\bq}\bq_l,V\otimes_{\bq}\bq_l)$ is a weak Mumford-Tate pair over $\bq_l$ ([Pin98]). If $d>3$, then $\MT_A\otimes_{\bq}\bq_l$ is of type $A$ and $V\otimes_{\bq}\bq_l$ is \textit{not} the standard representation. So the pair corresponds to column 4 of [Pin98] and hence, the Mumford-Tate conjecture holds for $A$.

If $d=3$ and $r=1$, we have $\rank \MT_A=3$. On the other hand, ([Orr15], Corollary 3.2) implies $18 = \dim V_l\leq 6\cdot 2^{\rank G_l-1} $ and hence, $\rank G_l\geq 3$. So $\rank G_l\geq \rank \MT_A=3$ and hence, the Mumford-Tate conjecture holds for $A$.\end{prf}

An important consequence of this is that away from a density zero set of places, the reduction $A_v$ is isogenous to the $d$-th power of a simple abelian variety. This follows from the following theorem of Zywina.

\begin{Thm} $\mr{([Zyw14],\; Theorem 1.2)}$ Let $A$ be an absolutely simple abelian variety over a number field $F$ such that the $F=F_A^{\conn}$ and the Mumford-Tate conjecture holds for $A$. Let $d^2$ be the dimension of the division algebra $\End^0(A)$ over its center. Then $A_v$ is the $d$-th power of an absolutely simple abelian variety away from a set of places of Dirichlet density zero.
\end{Thm}

\begin{Prop} Let $A$ be a simple abelian variety over a number field $F$ such that the $l$-adic monodromy groups are connected. Let $\End^0(A)$ be of dimension $d^2$ over its center $E$. Then, for all places $v$ in a set of density one, the reduction $A_v$ is isogenous to the $d$-th power $B_v^d$ of an abelian variety $B_v$ such that $\End^0(B_v)$ is a product of matrix algebras over CM fields.\end{Prop}

\begin{prf} Let $\mc{S}$ be the set of places $v$ of $F$ such that:\\
- $A_v$ has no supersingular abelian subvariety.\\
- $v$ has local degree one over $\bq$.

The first condition holds for a density one set of primes by ([Pin04], Proposition 1.6) and the second is a consequence of the Chebotarev density theorem. Thus, $\mc{S}$ has density one.

Let $\pi$ be a zero of $\PAv$ and let $B_{\pi}$ be the corresponding simple abelian variety. Since $B_{\pi}$ is not supersingular, it follows that the center $\bq(\pi)$ of $\End^0(B_{\pi})$ is a CM field. Furthermore, since $k_v$ is the prime field, $\End^0(B_{\pi})$ is commutative. Thus, $\End^0(B_{\pi}) = \bq(\pi)$. \end{prf}

\vspace{-0.3cm}
\begin{center}
\section{\fontsize{11}{11}\selectfont  Some preliminary results}
\end{center}
\vspace{-0.1cm}
 
In this section, we state and prove a few basic propositions that will be necessary in the subsequent sections. We begin with the following two facts about central simple algebras.

\begin{Thm}$\mr{([Jac80]\; Theorem\; 4.7)}$ Let $E$ be a local or global field. Let $A\hra B$ be an embedding of central simple algebras over $E$. Then $B\cong A\otimes_E \Cent_B(A)$ where $\Cent_B(A)$ is the centralizer of $A$ in $B$. Furthermore, $\Cent_B(A)$ is central simple over $E$.\end{Thm}

\begin{Corr} Let $E$ be a local or global field. Let $A\hra B$ be an embedding of central simple algebras over $E$. Then the integers $[A:E]$ and $\frac{[B:E]}{[A:E]}$ are relatively prime.\end{Corr}

\begin{prf} Since we have an embedding $A\hra B$, $C:=\Cent_B(A)$ is a central simple algebra over $E$ of dimension $\frac{[B:K]}{[A:K]}$. By the double centralizer theorem, $\Cent_B(C) = A$. 

Suppose, by way of contradiction, that there exists a prime $p$ dividing $[A:K]$ and $[C:K]$. Then there exists a field extension $E_1/E$ of degree $p$ such that $E_1$ can be embedded in both $A$ and $C$ (see the next proposition). Since any element of $E_1\bs E$ does not commute with all of $C$, we have a contradiction.\end{prf}

\begin{Prop} Let $K$ be a global or local field, $p$ a rational prime and let $C_1,C_2$ be central simple algebras over $K$ such that $p$ divides $[C_1:K]$ and $[C_2:K]$. Then there exists a field extension $L/K$ of degree $p$ that can be embedded in both $C_1$ and $C_2$.\end{Prop}
\begin{prf} We write $C_i = \Mat_{n_i}(D_i)$, $i=1,2$ where the $D_i$ are central division algebras over $K$. If $p|n_1$, then any extension of degree $p$ may be embedded in $C_1$. So any degree $p$ extension $L/K$ that has an embedding in $C_2$ satisfies the condition. So we may assume without loss of generality that $p\nmid n_1n_2$ and hence, $d_1^2:=[D_1:K]$, $d_2^2:=[D_2:K]$ are divisible by $p^2$. 

Let $\mc{P}$ be the set of primes of $K$ that either $D_1$ or $D_2$ is ramified at. Choose a cyclic extension $\wti{L}/K$ of degree $d_1d_2$ such that every prime in $\mc{P}$ is inert in $\wti{L}/K$. Let $L_1$, $L_2$, $L$ be the intermediate fields such that $[L_i:K] = d_i$ ($i=1,2$) and $[L:K]=p$. Then $L\sub L_i$ ($i=1,2$). Furthermore, $L_i$ splits $D_i$ are hence, is a maximal subfield of $D_i$.\end{prf}

\begin{Prop} An abelian variety $A$ of type $\IV(1,d,\wedge^r(\mr{Std}))$ has potential good reduction everywhere.\end{Prop}

\begin{prf} Since the derived group $\MT_A^{\der}$ is anisotropic over $\bq$, this follows from ([Per12], Thm 4.1.9).\end{prf}

\vspace{-0.2cm}
\begin{center}
\section{\fontsize{11}{11}\selectfont The density of ordinary reduction}
\end{center}

In this section, we show that the abelian varieties of type $\IV(1,d,\wedge^r(\mr{Std}))$ class have ordinary reduction at most places. We will draw heavily from the proof of ([Pin98], Theorem 7.1) which deals with the case where $\End(A) = \bz$ and the $l$-adic monodromy groups are of type $A$.

\subsection{\fontsize{11}{11}\selectfont Cocharacters}

For an abelian variety $A$ over a number field $F = F_A^{\conn}$, let $\MT_A$ be the Mumford-Tate group. Let $v$ be a place such that $\rank \Phi_{A_v} = \rank A$. We choose an element $t_v\in \GL_{2g}(\bq)$ with characteristic polynomial $\PAv$, which we defined as the characteristic polynomial of $\rho_l(\Fr_v)$ for any $l\neq p:=\charac(v)$. We denote by $\bT_v$ the Zariski closure of $t_v$ in $\GL_{2g}(\bq)$. 

The cocharacter group $X_*(\bT_v):=\Hom(\bg_{m,\qbar},\bT_{v}\times_{\bq} \qbar)$ and the character group $X^*(\bT_v):=\Hom(T_{v}\times_{\bq} \qbar,\bg_{m,\qbar})$ are in canonical perfect duality \vspace{-0.2cm} $$\la\;,\ra:X^*(\bT_v)\times X_*(\bT_v) \lra \bz,\;\;\;(\chi,\lamb)\mapsto \deg(\chi\circ\lamb).$$ When tensored with $\br$, this yields a perfect $\br$-valued pairing between the character space $X^*(\bT_v)\otimes_{\bz} \br$ and the cocharacter space $X_*(\bT_v)\otimes_{\bz} \br$.

We note that unlike the $p$-adic monodromy group $G_p(A)$, the subgroup $H_v:= \ov{\rho_{_p}(D_v)}^{\mr{Zar}}$ might not be reductive. In fact, in the case of ordinary reduction - which we expect to occur at a density one set of places - the group $H_v$ is solvable  and hence, cannot be reductive unless it is a torus. However, if $v$ is a place such that $\rank \Phi_{A_v} = \rank A$, then there exists an element $g_v\in \GL_{2g}(\qbar_p)$ such that \vspace{-0.15cm}$$g_v \bT_{v,\qbar_p}g_v^{-1} \sub H_{v,\qbar_p}\sub G_{p,\qbar_p}\;\;\;\;\mr{([Pin98], 3.12)}.$$ In particular, $H_v$ is of the same rank as $G_p(A)$ for such places.

\begin{Def}\normalfont A \textit{strong} Hodge cocharacter of $\bT_v$ is a cocharacter $\mu$ of $\bT_v$ such that there exists $g_v\in \GL_{2g}(\qbar_p)$ such that $g_v\mu$ is a Hodge cocharacter of $H_v$. A weak Hodge cocharacter is a cocharacter of $\bT_v$ that is conjugate under $\GL_{2g}(\qbar)$ to a strong Hodge cocharacter. \end{Def}

For a prime $\mfp$ of $K$ of local degree one over a $\bq$, let
$\ord_{\mfp}:K_{\mfp}\cong\qbar_p\lra \bq\cup \{\infty\}$ be the normalized valuation with $\ord_{\mfp}(\mfp)=1$. A \textit{quasi-cocharacter} of $\bT_v$ is a homomorphism $\mu:\bg_{m,\qbar}\lra \bT_{v,\qbar}$ such that $\mu^N$ is a cocharacter for some $N\in \bz$. In other words, the set of quasi-cocharacters of $\bT_v$ is $X_*(\bT_v)\otimes_{\bz} \bq$.

\begin{Def}\normalfont The Newton (quasi-)cocharacter $\eta_v$ is the unique quasi-cocharacter of $\bT_v$ such that for every $\chi\in X^*(\bT_v)$, $\la\chi,\eta_v\ra = \frac{\ord_p(\chi(t_v))}{[\bF _v:\bF _p]}$, where $\bF _v$ is the residue field at $v$.\end{Def}

By [LP97], the Frobenius tori $\bT_v$ for places $v$ with $\rank \Phi_{A_v} = \rank A$ have the same rank and the same formal character. Hence, we may fix a \textit{split} torus $\bT_0$ of $\GL_{2g,\bq}$ and conjugate it into every $G_{l,\qbar_l}$. For every rational prime $l$, we choose and fix an element $f_l\in \GL_{2g}(\qbar_l)$ such that \vspace{-0.15cm}$$\bT_{0,\qbar_l}\sub f_l^{-1}G_{l,\qbar_l}f_l.$$ We may assume that the torus $f_l^{-1}\bT_{0,\qbar_l}f_l$ is defined over $K_{\mf{l}}$ for some prime $\mf{l}$ of $K$ lying over $\bq_l$. In particular, if $l$ splits (completely) in $K$, we may assume $f_l^{-1}\bT_{0,\qbar_l}f_l$ is defined over $\bq_l$. Furthermore, for every place $v$ such that $\rank \Phi_{A_v} = \rank A$, we choose an element $h_v\in \GL_{2g}(\qbar)$ such that \vspace{-0.15cm}$$h_v^{-1}\bT_{v,\qbar}h_v = \bT_{0,\qbar}.$$

There exists a cocharacter $\eta_{_{0,v}}$ of $\bT_{0}$ that is conjugate to the Newton cocharacter $\eta_v$. Note that for every weight $\chi\in X^*(\bT_{0})$ that occurs in the representation $V(A)$, we have $\la \chi,\eta_{_{0,v}}\ra \leq 1$. 

\begin{Def} \normalfont For a quasi-cocharacter $\mu$ of $H_v$, we let $S_{\mu}$ denote the set of all\\ $H_v(\qbar_p)\rtimes \absqp$-conjugates of $\mu$ that factor through $\bT_{0,\qbar_p}$. \end{Def}

The following theorem relates ordinary reduction to these orbits of the Hodge and Newton cocharacters coinciding.

\begin{Thm} $\mr{([Pin98],\;Theorem\; 1.5)}$ $A$ has ordinary reduction at $v$ if and only if $S({\muH}) = S({\eta_v})$. \end{Thm}


\begin{Thm} $\mr{([Pin98]\; Theorem\; 3.15)}$ If $F=F_A^{\conn}$ and $v$ is a place of local degree one over $\bq$ such that $\rank \Phi_{A_v} = \rank A$, then $\eta_v$ lies in the convex closure of the orbit of $\mu_{0,v}$ under the Weyl group.\end{Thm}

\vspace{-0.5cm}
\subsection{\fontsize{11}{11}\selectfont Mumford-Tate pairs}
\vspace{-0.2cm}
Consider a connected reductive group $G$ over a field $E$ of characteristic zero and faithful finite dimensional representation $\rho$.

\begin{Def}\normalfont The pair $(G,\rho)$ is a \textit{weak} Mumford-Tate pair if there exist cocharacters $\mu_i:\bg_{m,\Ebar}\rightarrow G_{\Ebar}$ such that:\\
i. $G_{\Ebar}$ is generated by the images of all $G(\Ebar)$-conjugates of the $\mu_i$ and\\ 
ii. the weights of each $\rho\circ\mu_i$ are in $\{0,1\}$.\\
The pair is a \textit{strong} Mumford-Tate pair if these conditions hold and the cocharacters $\mu_i$ are conjugate under $\Gal_E$.\end{Def}

\noindent We list a few examples that are particularly important for our purposes in this article.

\noindent 1. The pair $(\MT_A,V(A))$ is a weak Mumford-Tate pair. This follows from the very definition.

\noindent 2. The pair $(G_l^0(A), V_l(A))$ for an abelian variety $A$ is a weak Mumford-Tate pair. This was shown in ([Pin98], Theorem 5.10).

\noindent 3. A pair $(G,\rho)$ is a weak Mumford-Tate pair if and only if $(\rho'(G),\rho')$ is one, for every irreducible summand $\rho'$ of $\rho$.

Let $A$ be an abelian variety of type $\IV(1,d,\wedge^r(\mr{Std}))$ over a number field $F$ with $D=\End^0(A)$ a central division algebra over $K$. Passing to a finite extension if necessary, we may assume $F = F_A^{\conn}$. The representation $V(A)_{\qbar}$ is a sum of $2d$ irreducible representations. Let $V_1$ be one of these irreducible components and set $\rho_1:=\rho|_{V_1}$. Then $(\MT_A,\rho_1)$ is a weak irreducible Mumford-Tate pair.

\subsection{\fontsize{11}{11}\selectfont The groups $W(G,T)$ and $\Pi(G,T)$:} Let $G$ be a reductive group over a perfect field $E$. Its maximal tori are all conjugate to each other over the algebraic closure $\ov{E}$. Choose a maximal torus $T$ of $G$. We have a homomorphism \vspace{-0.1cm} $$\phi_{_T}:\abse\lra \Aut(X^*(T)),\;\;\;\;\; \si\mapsto\; (\chi \mapsto \si(\chi))$$ for $\si\in \abse$, $\chi\in X^*(T)$. Note that $T$ is split if and only if $\phi_{_T}(\abse) = \{1\}$.

We use the symbol $W(G,T)$ to mean the Weyl group $N_{G}(T)(\ov{E})/T(\ov{E})$ where $N_{G}(T)$ is the normalizer of $T$ in $G$. Now, $\absk$ acts on $W(G,T)$ and $\si(w) = \phi_{_T}(\si)\circ w\circ \phi_{_T}(\si)^{-1}$ for $\si\in \absk$, $w\in W(G,T)$. We define $\Pi(G,T)$ to be the subgroup of $\Aut(X^*(T))$ generated by $W(G,T)$ and $\phi_{_T}(\abse)$. As abstract groups, $W(G,T)$ and $\Pi(G,T)$ are independent of $T$. We note that while $W(G,T)$ is insensitive to base change, the groups $\phi_{_T}(\abse)$ and $\Pi(G,T)$ are not.

Suppose we have a faithful representation $\rho:G\hra \GL(V)$ for some finite dimensional vectors space $V$ over $E$. We denote by $\Om(V)\sub X^*(T)$ the set of weights of $V$. We have a decomposition \vspace{-0.15cm}$$ V\otimes_E \ov{E} = \bigoplus\limits_{\chi\in \Om(V)} V(\chi)$$ and hence, for every $t\in T(\ov{E})$, we have \vspace{-0.15cm}$$\det(XI - \rho(t)) = \prod\limits_{\chi\in \Om(V)} (X-\chi(t))^{m_{\chi}}$$ where $m_{\chi}:=\dim V(\chi)$ is the multiplicity of $\chi$.

We are primarily concerned with the case where the reductive group $G$ is the Mumford-Tate group of an abelian variety $A$ of type $\IV(1,d,\wedge^r(\mr{Std}))$ with $D:=\End^0(A)$ a central division algebra over an imaginary quadratic field $K$. Let $\wti{D}$ be the unique central division algebra over $K$ such that $[D]^{-1} = [\wti{D}]^r$ in $\Br(K)$. Fix a maximal torus $\bT$ in $\MT_A$. Since $\MT_A$ is a $\bq$-form of $\GL_d$, the Weyl group $W(\MT_A,\bT)$ is isomorphic to the symmetric group $S_d$. Furthermore, a field extension $L/K$ of degree $d$ splits $\bT$ if and only if $L$ is a CM field and splits $\wti{D}$. As seen in the next proposition, this is equivalent to $L$ splitting $D$. 

\begin{Prop} Let $K$ be a number field and $D$ a $d^2$-dimensional central division algebra over $K$. Let $r$ be an integer prime to $d$ and let $\wti{D}$ be the unique central division algebra over $K$ such that $[\wti{D}] = [D]^r$ in $\Br(K)$. Then any field extension $L/K$ has an embedding in $D$ if and only if it has an embedding in $\wti{D}$.\end{Prop}

\begin{prf} Any subfield of $D$ or $\wti{D}$ lies in a maximal subfield. So we may assume without loss of generality that $[L:K] = d$.

Let $\mc{P}$,$\wti{\mc{P}}$ be the sets of primes of $K$ that $D$, $\wti{D}$ are ramified at. Note that for any $\mfp\in \mc{P}$, $\inv_{\mfp}(\wti{D}) = r\inv_{\mfp}(D)$ and since $r$ is prime to $d$, which is the order of $[D]$ in $\Br(K)$, it follows that $\mc{P} = \wti{\mc{P}}$. Now, any extension $L/K$ of degree $d$ has an embedding in $D$ (resp. $\wti{D}$) if and only if $L$ splits $D$ (resp. $\wti{D}$). But this happens if and only if for every prime $\mfp'$ of $L$ lying over $\mfp$,  the local degree $[L_{\mfp'}:K_{\mfp}]$ annihilates $\inv_{\mfp}(D)$ (resp. $\inv_{\mfp}(\wti{D})$). Since $\mc{P} = \wti{\mc{P}}$, the proposition follows.\end{prf}

\begin{Prop} Let $d$ be an odd integer. Let $D$ be a $d^2$-dimensional central division algebra over an imaginary quadratic field $K$. There exists infinitely many cyclic extensions $E/\bq$ of degree $d$ such that:\\
\noindent $\mr{i}$. $E$ is linearly disjoint from $K$\\
\noindent $\mr{ii}$. $EK$ is a CM field.\\
\noindent $\mr{iii}$. $EK$ has an embedding in $D$\end{Prop}

\begin{prf} Note that any cyclic extension of odd degree over $\bq$ is totally real and hence, its compositum with $K$ is a CM field cyclic over $\bq$.

Let $\mfp_1,\cdots,\mfp_r$ be the primes of $K$ that $D$ is ramified at and let $p_1,\cdots,p_r$ be the rational primes they lie over. There are infinitely many cyclic extensions $E/\bq$ of degree $d$ such that $p_1,\cdots,p_r$ are inert in $E/\bq$ [Th16, Lemma 6.2]. Since these extensions are Galois over $\bq$ and $d$ is an odd integer, any such extension $E$ is linearly disjoint from $K$. This ensures that $EK/K$ is cyclic of degree $d$. Furthermore, the primes $\mfp_1,\cdots,\mfp_r$ are inert in $EK/K$ since $p_1,\cdots,p_r$ are inert in $E/\bq$. So $[D]\in \Br(EK/K)$ and hence, $EK$ has an embedding in $D$.\end{prf}

Since there exists a cyclic extension $L/\bq$ of degree $2d$ that contains $K$ and splits the division algebra $\wti{D}$ (by Proposition 3.5), it follows that $\phi_{_{\bT}}(\absq)$ is cyclic of order $2d$ and\\ $\phi_{_{\bT}}(\absq)\cap W(\MT_A,\bT)$ is cyclic of order $d$. Thus, $\Pi(\MT_A,\bT)\cong S_d\times (\bz/2\bz)$. As noted earlier, these abelian varieties fulfill the Mumford-Tate conjecture. So [Zyw14] implies that away from a density zero set of places, $\GalW = \Pi(\MT_A,\bT)$. The fixed field of the Weyl group is $K$.

\begin{Prop} Let $A$ be an abelian variety of type $\IV(1,d,\wedge^r(\mr{Std}))$ over a number field $F$ such that $F = F_A^{\conn}$. At most places $v$ of good reduction, $\End^0(A_v) = \Mat_{d}(E_v)$ for some CM field $E_v$ of degree $2 {d\choose r}$ containing $K$.\end{Prop}

\begin{prf}	Since $A$ satisfies the Mumford-Tate conjecture, it is immediate from the main theorem of [Zyw14] that away from a set of density zero, $\GalW = \Pi(\MT_A)$. For any such place, the fixed field of $W(\MT_A)$ is $K$. Furthermore, by Proposition 1.9, there exists a density one set of places such that $\End^0(A_v)$ is a product of matrix algebras over CM fields. This completes the proof.\end{prf}

Since the pair $(\MT_A,V(A))$ is a weak Mumford-Tate pair, it follows that for any irreducible $\MT_{A,\qbar}$-representation $V_i\sub V(A)\otimes_{\bq} \qbar$, the pair $(\MT_{A,\qbar},V_i)$ is an irreducible weak Mumford-Tate pair. The Weyl group $W(\MT_A,\bT)$ acts transitively on the set $\Om_i$ of weights of $V_i$ and $\Pi(\MT_A,\bT)$ acts transitively on the disjoint union $\Om = \Om_1\cup \Om_2$. So each $\Om_i$ is stable under the action of $\absk$. 

Choose a prime $l$ that splits (completely) in $K/\bq$. Let $\mfl_1$, $\mfl_2$ be the primes of $K$ lying over $l$. Then $K\otimes_{\bq}\bq_l = K_{\mfl_1}\times K_{\mfl_2}$ and if we set $V_{\mfl_i}(A):=V_l(A)\otimes_{K\otimes\bq_l} K_{\mfl_i}$, we have a decomposition \vspace{-0.15cm}$$V_l(A) =  V_{\mfl_1}(A)\oplus V_{\mfl_2}(A)$$ of $\bq_l[\absf]$-modules. 

\vspace{-0.15cm}
\subsection{\fontsize{11}{11}\selectfont The variety $\Conj'(G)$}
\vspace{-0.1cm} 
Let $G$ be a connected reductive group over a local or global field $E$. The derived group $G^{\der}_{\Ebar}$ is an almost direct product of almost simple groups. We assume these simple groups are type $A$,$B$,$C$ or $D^{\mb{H}}$. 

We define $\mc{A}$ to be the group of automorphisms $\tau$ of the group $G_{\Ebar}$ such that: 

\noindent $\bullet\;$ $\tau$ is the identity on the center of $G_{\Ebar}$.\\
$\bullet\;$ Every simple factor is stable under $\tau$.\\
$\bullet\;$ The restriction of $\tau$ to any simple factor $H$ agrees with conjugation by some element of $H$.

Let $R_G$ be the affine coordinate ring of $G$. This is an $E$-algebra and the group $\mc{A}$ acts on $R_G$ by composition. We define $R_G^{\mc{A}}$ to be the $E$-subalgebra of elements fixed by $\mc{A}$. Furthermore, we define $\Conj'(G):=\mr{Spec}(R_G^{\mc{A}})$. This is a universal categorical quotient and is a $E$-form of the quotient $T_{\Ebar}/W(G,T)$ for any maximal torus $T$.

In this article, we are primarily concerned with the case where $G$ is either the Mumford-Tate group $\MT_A$ of an abelian variety or the connected component $G_l^0(A)$ of the $l$-adic monodromy group.
\vspace{-0.15cm}
\subsection{\fontsize{11}{11}\selectfont The set $\mc{S}_A$} 
\vspace{-0.15cm}
We choose a rational prime $l$ that that splits in $K/\bq$ and is prime to the discriminant of $D$ (and hence, prime to the discriminant of $\wti{D}$). We have the decomposition \vspace{-0.15cm}$$V_l(A) \cong V_{\mfl_1}(A)\oplus V_{\mfl_2}(A)$$ of $\bq_l[\absk]$-modules. Hence, the set $\Om(V_l(A))$ of weights of $V_l(A)$ is the union $\Om_1 \cup \Om_2$ where $\Om_i$ is the set of weights on $V_{\mfl_i}(A)$ . Furthermore, the Weyl group acts transitively on each $\Om_i$ and $\Pi(G_l)$ acts transitively on $\Om$. 

For every $\chi\in \Om_1$, we have \vspace{-0.15cm}$0\leq \la \chi,\eta_{_{0,v}} \ra \leq 1$. For all places $v\in \Sigma_F$ such that $\rank \Phi_{A_v} = \rank A$, we define $$a_{v,j}: = \mr{Tr}\left(\rho_{\mfl_j}(h_v^{-1}t_vh_v)\right)\cdot \mr{Tr}\left(\rho_{\mfl_j}(h_v^{-1}t_vh_v)^{-1}\right).$$ We first show that this element lies in $K$.

\begin{Prop} The algebraic numbers $a_{v,j}$ lie in $K$.\end{Prop}

\begin{prf} The elements $\chi(h_v^{-1}t_vh_v)$ lie in the multiplicative group $\Phi_{A_v}$ generated by $\WA$ and hence, the sum $a_{v,j}$ lies in $\bqW$. Since the set $\Om_j$ is stable under the action of the Weyl group, it follows that $a_{v,j}$ lies in the subfield of $\bqW$ fixed by the Weyl group $W(\MT_A)\cong S_d$. But as noted earlier, the fixed field of the Weyl group is $K$. Hence, $a_{v,j}\in K$. \end{prf}

\begin{Corr} The image of the map \vspace{-0.15cm}$$\tau:\bT_0\lra \mb{A}_{\qbar}^1,\;\;\; t\mapsto \sum\limits_{\chi\in \Om_j} \chi(t)$$ lies in $\mb{A}_{K}^1$.\end{Corr}

\begin{prf} We have $\tau(t_v)\in K$ for all $\{v:\rank \Phi_{A_v} = \rank A \}$. Since the elements $\{t_v:\rank \Phi_{A_v} = \rank A \}$ generate a dense subgroup of $\bT_0$, this corollary follows.\end{prf}

For an abelian variety $A$ of type $\IV(1,d,\wedge^r(\mr{Std}))$, we define the set $\mc{S}_A$ as the set of places $v$ of $F$ such that:\\
		1. $v$ has local degree one over $\bq$.\\
		2. The group $\Phi_{A_v}$ generated by $\WA$ is of full rank.\\
		3. $\Nm_{K/\bq}{a_v}\notin p\bz$.\\
		4. $\GalW = \Pi(\MT_A)$.
		
The general idea is to show that $\SA$ has Dirichlet density one and that each place $v\in \SA$ is a place of ordinary reduction. To this end, we will need the following two lemmas.
		
\begin{Lem} For any constant $C$, the set $\{\be\in \mc{O}_K:\Nm_{K/\bq}(\al)\leq C \}$ is finite.\end{Lem}
\begin{prf} Clearly, for any $\be\in \mc{O}_K$, the norm $\Nm_{K/\bq}(\be)\in \bz$. So suffices to show that the set $\{\be\in \mc{O}_K:\Nm_{K/\bq}(\be)=N \}$ is finite for any integer $N	\leq C$.

For any integer $N\leq C$, if $N =\Nm_{K/\bq}(\be)$, then $\be$ generates some ideal $I$ in $\mc{O}_K$ dividing $(N)$. Clearly, there are finitely many such ideals. For any such ideal $I$, if $I = \be\mc{O}_K = \be_1\mc{O}_K$, then $\frac{\be_1}{\be}$ is a unit in $\mc{O}_K$. But since $K$ is an imaginary quadratic field, any unit of $\mc{O}_K$ is a root of unity of order dividing $6$. This complete the proof.\end{prf}

\begin{Lem} For any $v\in \SA$, there exist weights $\chi_{_1},\chi_{_2}\in \Om_j$ such that $\la \chi_{_1}\chi_{_2}^{-1},\eta_{_{0,v}} \ra =-1$.\end{Lem}

\begin{prf} We show this for $\Om_1$. The proof for $\Om_2$ is identical. 

First, we note that $K\sub F$ since all the endomorphisms of $A$ are defined over $F_A^{\conn}$. Since $p:=\charac(v)$ splits completely in $F$, it is clear that $p$ splits in $K/\bq$. Let $\mfp$ be a prime of $K$ lying over $p$ and let $\ord_{\mfp}$ be a valuation extending $\ord_p$.

Let $\chi\in \Om_1$ be any weight that occurs in the absolutely irreducible representation $V_1$. The algebraic number $\chi(h_v^{-1}t_vh_v)$ is a unit at all primes not above $p$ and the archimedean norms are all equal to $\sqrt{p}$. Hence, any quotient $\chi_{_1}(h_v^{-1}t_vh_v)\chi_{_2}(h_v^{-1}t_vh_v)^{-1}$ is a unit away from $p$ and the archimedean norms are all equal to $1$. Hence, $a_{v,1}$ is an element of $K$ with norm $\Nm_{K/\bq}(a_{v,1})\leq p\big|\Om_1\big|^2$ and is integral away from $p$. Since $a_{v,1}\notin \mc{O}_K$, it follows that $\ord_p(\Nm_{K/\bq}(a_{v,1}))\leq -1$. So there exists a prime $\mfp$ of $K$ lying over $p$ such that $\ord_{\mfp}(\Nm_{K/\bq}(a_{v,1}))\leq -1$ and hence, either $\ord_{\mfp}(a_{v,1})\leq -1$ or $\ord_{\mfp}(\ov{a_{v,1}})\leq -1$.
 
We may assume without loss of generality that $\ord_{\mfp}(a_{v,1})\leq -1$. Hence, there exist two weights $\chi_{_1},\chi_{_2}\in \Om_1$ such that \vspace{-0.15cm}$$\la \chi_{_1}\chi_{_2}^{-1},\eta_{_{0,v}}\ra = \ord_{\mfp}\left(\chi_{_1}\chi_{_2}^{-1}\right)\leq -1.$$ On the other hand, we have \vspace{-0.15cm}$$\la \chi_{_1}\chi_{_2}^{-1},\eta_{_{0,v}} \ra = \la \chi_{_1},\eta_{_{0,v}}\ra-\la\chi_{_2},\eta_{_{0,v}}\ra\geq 0 - 1 = -1,$$ which forces equality.\end{prf}

\begin{Prop} Let $A$ be an abelian variety of type $\IV(1,d,\wedge^r(\mr{Std}))$ over a number field $F$ such that $F = F_A^{\conn}$. Then the set $\SA$ has Dirichlet density one.\end{Prop}

\begin{prf} The set of places of $F$ with local degree one over $\bq$ has density one by the Chebotarev density theorem. By [LP97], the set of places such that the inequality $\rank \Phi_{A_v} = \rank A = d$ has density one. So it suffices to show that the set \vspace{-0.15cm}$$\{v\in \Sigma_F :\Nm_{K/\bq}(a_{v,j})\notin p\bz\}$$ has density one.

Fix a rational prime $l$ that splits in $K/\bq$. Let $\mf{l}$ be a place of $K$ lying over $l$. Then $K_{\mf{l}}\cong \bq_l$. As before, let $\Conj'(f_l^{-1}G_lf_l)$ be the affine $\bq_l$-variety of semisimple conjugacy classes, which is a $\bq_l$-form of the quotient $f_l\bT_{0,\qbar_l}f_l^{-1}/W(G_l,f_l\bT_{0,\bq_l}f_l^{-1})$. We write $W_l:=W(G_l,f_l\bT_{0}f_l^{-1})$ for brevity.

We define \vspace{-0.15cm}$$G_l\lra \bT_{0,\bq_l}/W_l,\;\;\; g\mapsto f_l^{-1}g^{\mr{ss}}f_l \pmod {W_l} $$ where $g^{\mr{ss}}$ is the semisimplification of $g$. This yields a dominant morphism \vspace{-0.15cm}$$G_{l,\qbar_l}\lra \Conj'(f_l^{-1}G_lf_l)_{\qbar_l}$$ of affine varieties when base changed to the algebraic closure $\qbar_l$. Furthermore, the morphism \vspace{-0.15cm}$$\bT_0\lra \mb{A}_K^1,\;\;\;t\mapsto \mr{Tr}(t)\cdot \mr{Tr}(t^{-1})$$ factors through $\bT_0/W$ since $\Om_1$ is stable under the action of the Weyl group. So we may form the compositions $$\psi_j:G_l\lra \mb{A}_{K_{\mf{l_j}}}^1\cong \mb{A}_{\bq_l}^1,\;\;\;\;j=1,2   .$$ 

Since the representation $V_j$ is absolutely irreducible, the linear independence of characters implies that $\psi_j$ is non-constant. We define the sets\vspace{-0.15cm}$$X_{\mfl_j}:=\{g\in G_l:\psi(g)\in \mc{O}_K\text{ and } \Nm_{K/\bq}\left(\psi(g)\right)\leq |\Om_j|^4\}.$$ and $X_l:=X_{l,1}\cup X_{l,2}$. As seen in the preceding lemma, this allows for only finitely many values for $\psi(g)$ and hence, $X_l$ is a nowhere dense Zariski subset of $G_l$. Now, let $v$ be a place of $F$ such that:\\
		-$v\nmid l$\\
		-$v$ has local degree one over $\bq$\\
		-$\rank \Phi_{A_v} = \rank A$

Then, by construction, $a_{v,j}=\psi_j(\rho_{\mfl_j}(t_v))$ and $v\in \mc{S}_A$ if and only if $\rho_{\mfl_j}(t_v)\notin X_l$. But $\rho_{\mfl_j}(\absf)$ is a compact $l$-adic analytic subgroup of $G_l(\bq_l)$ which, by definition, is Zariski dense in $G_l(\bq_l)$. Hence, $Y_{\mfl_j}:=X_{\mfl_j}(\bq_l)\cap \rho_{\mfl_j}(\absf)$ is a nowhere dense closed analytic subset of $\rho_{\mfl_j}(\absf)$.

Let $\mu_j$ be the Haar measure on $\rho_{\mfl_j}(\absf)$ with total volume one. Then $\mr{vol}_{\mu}(U)\lra 0$ as $U$ runs through a cofinal system of open compact neighborhoods of $Y_{\mfl_j}$. From the Chebotarev density theorem, it follows that each $Y_{\mfl_j}$ and hence, $Y_l:=Y_{\mfl_1}\cup Y_{\mfl_2}$ have volume zero. This completes the proof.\end{prf}

\begin{Thm} Let $A$ be an abelian variety of type $\IV(1,d,\wedge^r(\mr{Std}))$ over a number field $F$ such that $F = F_A^{\conn}$. Then $A$ has ordinary reduction at a set of places of density one.\end{Thm}

\begin{prf} Since $\SA$ has Dirichlet density one, it suffices to show that every place $v\in \SA$ is a place of ordinary reduction. Let $v\in \mc{S}_A$ and let $\eta_{_{0,v}}$ be the Newton cocharacter. We choose an element $t_v$ with characteristic polynomial $\PAv$ and let $\bT_v$ be its Zariski closure in $\GL_{2g}(\bq)$.

We take a strong Hodge cocharacter of $\bT_v$ and conjugate it via $h_v$ into a cocharacter $\mu_{0,v}$ of $\bT_v$. By ([Pin98], Theorem 3.15), the cocharacter $\eta_{_{0,v}}$ lies in the convex closure of the orbit of $\mu_{0,v}$ under the Weyl group. Let $X_*(\bT_0)\otimes_{\bz}\br =X_0\oplus X_1$ be a decomposition of the cocharacter space with $X_0$ belonging to the central part. Since $\eta_{_{0,v}}\in X_*(\bT_0)\otimes_{\bz}\br$, we may write $\eta_{_{0,v}} = (\eta_{_0},\eta_{_1})$ with $\eta_i\in X_i$, $i=0,1$. By the convex closure theorem, $\eta_{_0} = \mu_{_0}$ and $\eta_{_1}$ lies in the convex closure of the $W_l$-orbit of $\mu_1$.
	
We may identify the set $\Om_1$ of weights with a subset of $\br^d$ whose linear span is \vspace{-0.15cm}$$\{(x_1,\cdots,x_{d})\in \br^d:\sum\limits_{i=1}^{d} x_i = 0\}$$ and the characters $\chi_{_1},\chi_{_2}$ are mapped to $e_1,e_2$. So the character space $X^*(\bT_v)\otimes_{\bz} \br$ is the subspace \vspace{-0.15cm}$$\{(x_1,\cdots,x_{d})\in \br^d:\sum\limits_{i=1}^{d} x_i = 0\}$$ and the cocharacter space $X_*(\bT_v)\otimes_{\bz} \br$ may be identified with the quotient space $\br^d/\br(1,\cdots,1)$. We denote the elements of the cocharacter space by $[x_1,\cdots,x_d]/{\sim}$.
	
Now, we have $\mu_{_{0,h}} = [1,0,\cdots,0]/{\sim}\in X_*(\bT_v)$ up to conjugation by $W_l\cong S_d$. Furthermore, the cocharacter $\eta_{_{0,v}}$ lies in the convex closure of the Hodge cocharacter. So we may write $\eta_{_1} = [x_1,\cdots,x_d]/{\sim}$ where $1\geq x_1\geq\cdots\geq x_d \geq 0$ and $\sum\limits_{i=1}^{d}x_i = 1$. Then \vspace{-0.15cm}$$1 = \la e_1-e_2, \muN  \ra = y_1-y_2 \leq y_1\leq 1 $$ and hence, $y_1=1$, $y_2=0$.

Thus, the Newton and Hodge cocharacters lie in the same orbit of the Weyl group. Hence, $A_v$ is ordinary by ([Pin98], Theorem 1.5).\end{prf}

\begin{center}
	\section{\fontsize{11}{11}\selectfont Splitting of reductions of type $\IV(1,3,\mr{Std})$}
\end{center}

In this section, we study the possible splitting types of the reduction of an abelian variety of type $\IV(1,3,\mr{Std})$. Since these abelian varieties fulfill the Mumford-Tate conjecture, [Zyw14] implies that their reductions are third powers of simple three-folds away from a set of Dirichlet density zero. We study the types of splitting at these exceptional set of density zero which might be an infinite set.

The following is a well known fact about the reductions of simple abelian varieties. We provide a proof for the reader's convenience.

\begin{Prop} Let $A$ be a simple abelian variety over a number field $F$. Let $v$ be a place of good reduction and let $A_v =_{\text{isog}} \prod\limits_{i=1}^k B_i^{d_i}$ be a simple decomposition up to isogeny. Then we have an injection $\End^0(A)\hra \Mat_{d_i}(\End^0(B_i))$ for every $i$.\end{Prop}

\begin{prf} We have the map $A\lra A_v\xra{\mr{ projection }}B_i^{d_i}$ which induces a map $\phi_i:D\lra \Mat_{d_i}(\End^0(B_i)) $. Suppose this is not an injection and let $\al\in D$ be an element such that its image $\phi_i(\al) =0$ in $\Mat_{d_i}(\End^0(B_i))$. Choose an element $\hat{\al}$ such that $\hat{\al}\circ \al= [N]\in D$ is the multiplication by $N$ map for some integer $N\neq 0$. Then $\phi_i([N]) = 0$ and hence, $\phi_i([1]) =0$. Since $\phi_i([1]) = [1]\in \Mat_{d_i}(\End^0(B_i))$, this is a contradiction.\end{prf}

A consequence of [Zyw14] is that an abelian variety of type $\IV(1,3,\mr{Std})$ has reduction isogenous to the third power of a simple three-fold away from a set of density zero. But the reductions might exhibit different behavior at this density zero (but possibly infinite) set of places.

Note that $\absq$ acts transitively on $\Om$ and hence, acts on $\Om$ through a transitive subgroup of $S_d\times S_2$. In particular, if $d$ is a prime, the projection of the image to $S_d$ contains a $d$-cycle $\sigma_d$. By the Chebotarev density theorem, there exist infinitely many primes $l$ such that $\sigma_d\in \absql$. 

Consider the $\absq$-equivariant map \vspace{-0.15cm}$$\ev:X^*(\bT)\lra \Phi_{A_v},\;\;\;\chi\mapsto \chi(t_v).$$ Now, $X^*(\bT)\cong \bz^{d-1}$. We may assume without loss of generality that the $d$-cycle is $(1\;2\;\cdots\; d)$. The set of sublattices of $\bz^{d-1}$ invariant under this cycle is in bijection with the ideals of $\bz[X]/(\frac{X^d-1}{X-1})$. The ring $\bz[X]/(\frac{X^d-1}{X-1})$ is isomorphic to $\bz[\zeta_d]$ and hence, the only possibilities for $\ker(\ev)$ are the following:\\
-$0$ \\
-$\bz^{d-1}$\\
-$\{(x_1,\cdots,x_d):\sum\limits_{i=1}^d x_i=0 \}$

In the first case, the eigenvalues of the Frobenius are distinct. In the second case, all the eigenvalues are equal to $\sqrt{q}$ and hence, the reduction is supersingular. In the third case, there are precisely two distinct eigenvalues and hence, the reduction is of type $\IV(1,d)$.

For the rest of this section, let $A$ be an abelian variety of type $\IV(1,3,\mr{Std})$ over a number field $F$. Let $v$ be a place of good reduction such that $p=\charac(v)$ does not divide the discriminant of $D := \End^0(A)$. 

Replacing $F$ by a suitable finite extension if necessary, we may assume the following:\\
-the size $q$ of the residue field is a square.\\
-the $l$-adic monodromy groups are connected.\\
-the group $\Phi_{A_v}\sub \qbar^*$ is torsion-free.

\begin{Prop} $A_v$ is isogenous to the third power $B_v^3$ of an abelian variety $B_v$ over $\bF _v$.\end{Prop}

\begin{prf} We show that for any simple component $B$, the maximum power of $B$ in the simple decomposition of $A_v$ is a multiple of $3$. 

The polynomial $\PAv$ is the third power $f_{A_v}(X)^3$ for some polynomial $f_{A_v}\in \bz[X]$. Passing to a suitable finite extension if necessary, we may assume the Weil numbers corresponding to supersingular components are rational. So the largest integer $t$ such that $(X-q)^t|\PAv$ is a multiple of $3$ and hence, the supersingular part of $A_v$ has dimension divisible by $3$. 

Let $\pi\notin \bq$ be a zero of $f_{A_v}$ and let $B_{\pi}^t$ be the corresponding iso-simple abelian subvariety of $A_v$. Since $B_{\pi}$ is not supersingular, the field $\bq(\pi)$ is a CM field. Suppose, by way of contradiction, that $t\nmid 3$. Then $3$ divides $[D_{\pi}:\bq(\pi)]$. Thus, we have the following possibilities:\\
 (i) $[\bq(\pi):\bq]=2$, $[D_{\pi}:\bq(\pi)]=9$\\
(ii) $[\bq(\pi):\bq]=2$, $[D_{\pi}:\bq(\pi)]=81$\\
(iii) $[\bq(\pi):\bq]=6$, $[D_{\pi}:\bq(\pi)]=9$

Note that in all three cases, $K\hra D\hra D_{\pi}$ and if $K\nsub\bq(\pi)$, then $\bq(\pi)K$ is a degree $4$ or degree $12$ subfield of $D_{\pi}$, a contradiction since $4\nmid[D_{\pi}:\bq]$. Hence, $K\sub\bq(\pi)$.

In case (i), the dimensions of $D$ and $D_{\pi}$ coincide and hence, the embedding $D\hra D_{\pi}$ is an isomorphism. But this is not possible unless $D$ is ramified only at the places of $K$ lying over $p$. 

In case (ii), $D$ and $D_{\pi}$ are central over the same field $K = \bq(\pi)$ and hence, Corollary 2.2 yields a contradiction.

In case (iii),  $D\otimes _{K}\bq(\pi)\hra D_{\pi}$ and since the dimensions coincide, $D\otimes _{K}\bq(\pi)\cong D_{\pi}$. Since $D$ is split at the places of $K$ lying over $p$, it follows that $D_{\pi}$ is split at the places of $\bq(\pi)$ lying over $p$. This is a contradiction, since $[D_{\pi}]\neq 0$ in $\Br(\bq(\pi))$.\end{prf}

\begin{Prop} The reduction $A_v$ does not have an ordinary elliptic curve as an abelian subvariety. \end{Prop}

\begin{prf} Suppose $A_v$ has an ordinary elliptic curve $E$ as a simple component and let $t$ be the multiplicity with which it occurs. The endomorphism algebra $K_1:=\End^0(E)$ is an imaginary quadratic field. We have the injection $D\hra \Mat_{t}(K_1)$ and hence, $3$ divides $t$. Furthermore, $D$ and $\Mat_{t}(K_1)$ are of the same dimension and hence, $t=6$ or $9$.

Let $\pi_i,\ov{\pi}_i$ ($i=1,2,3$) be the zeros of $f_{A_v}$ and set $\al_i=\frac{\pi_i}{\sqrt{q}}$. We may order the $\pi_i$ so that $\al_1\al_2=\al_3$. Since one of the corresponding elliptic curves occurs with multiplicity $\geq 6$ in the simple decomposition of $A_v$, we may assume $\al_1=\al_2$. So $\al_1^2=\al_1$ and since $\Phi_{A_v}$ is torsion-free, it follows that $\al_1=\al_2=\al_3=1$. So $A_v$ is supersingular, a contradiction.\end{prf}

\begin{Prop} If $A_v$ has a supersingular simple component, then $A_v$ is supersingular.\end{Prop}

\begin{prf} Suppose $A_v$ has the supersingular elliptic curve $E$ as a simple component and let $t$ be the multiplicity with which it occurs. Then $t=3$ or $6$ and we have the injection $D\hra \Mat_t(\bq_{p,\infty})$ where $\bq_{p,\infty}$ is the unique quaternion algebra over $\bq$ ramified only at $p$ and $\infty$. 

We have the surjective homomorphism $X^*(\bT_l)\lra \WA$ that maps the set of weights on to $\WA$. Let $\pi_i,\ov{\pi}_i$ ($i=1,2,3$) be the zeros of $f_{A_v}$. We may order the $\pi_i$ so that $\pi_1\pi_2=\pi_3\sqrt{q}$. So, if $\pi_1 = \pi_2 = \sqrt{q}$, then $\pi_3 = \sqrt{q}$ as well. Hence, the supersingular part of $A_v$ is not of dimension $6$.

Suppose, by way of contradiction, that $t=3$. Then $D\hra \Mat_3(\bq_{p,\infty})$. In particular, we have the embedding $K\hra \Mat_3(\bq_{p,\infty})$ and since $2=[K:\bq]$ is prime to $3$, we have $K\hra \bq_{p,\infty}$. Hence, $K$ splits $\bq_{p,\infty}$. Now, $D\hra \Mat_3(\bq_{p,\infty})\otimes_{\bq} K\cong \Mat_{36}(K)$. Let $C$ be the centralizer of $D$ in $\Mat_{36}(K)$. Then $D\otimes_{K} C\cong\Mat_{36}(K) $ and hence, $[C:K] =4$. But since $D$ is $3$-torsion in $\Br(K)$, so is $C$, a contradiction.\end{prf}

The following corollaries now follow immediately.

\begin{Corr} Let $A$ be an abelian variety of type $\IV(1,3,\mr{Std})$ over a number field $F$. Let $v$ be a place of good reduction such that $p=\charac(v)$ does not divide the discriminant of $D := \End^0(A)$. Then $A_v$ is isogenous to one of the following:\\
$1$. $A_{1}^9$ where $A_{1}$ is the supersingular elliptic curve.\\
$2$. $B_v^3$ where $B_v$ is absolutely simple of dimension $3$.
\end{Corr}

\begin{Corr} The endomorphism algebra $\End^0(A_v)$ has the following possibilities:\\
$1$. $\Mat_3(L)$ where $L$ is a CM field of degree $6$.\\
$2$.  $\Mat_3(D_1)$ where $D_1$ is a $9$-dimensional central division algebra over an imaginary quadratic field\\
$3$. $\Mat_{9}(\bq_{p,\infty})$ where $\bq_{p,\infty}$ is the quaternion over $\bq$ ramified only at $p$ and $\infty$. \end{Corr}

\noindent Note that case 2 occurs precisely when $B_v$ is simple with Newton polygon $3\times\frac{1}{3},3\times\frac{2}{3}$.


\begin{center}
\section{\fontsize{11}{11}\selectfont Newton Polygons of type $\IV(1,3,\mr{Std})$}
\end{center}
\vspace{-0.3cm}

\noindent \textbf{Notation:} For a Weil $q$-integer $\pi$ (up to conjugacy) we denote the corresponding simple abelian variety (up to isogeny) over $\mb{F}_q$ by $B_{\pi}$. We denote the endomorphism algebra $\End(B_{\pi})$ by $D_{\pi}$. This is a division algebra central over $\bq(\pi)$ since $B_{\pi}$ is simple. Furthermore, the Hasse invariants of $D_{\pi}$ are given by \vspace{-0.2cm} $$\inv_w(D_{\pi}) = \begin{cases} \frac{1}{2} & \text{ if }w \text{ is archimedean}\\
0 & \text{ if }w\nmid p\\
\frac{w(\pi)}{w(q)} [\bq(\pi)_w:\bq_p]\;\pmod{1} & \text{ if } w|p.\end{cases}$$

The Newton slopes of $B_{\pi}$ are given by $\frac{w(\pi)}{w(q)}$ where $w$ runs through the places lying over $p$. As before, $\mc{W}_{B_{\pi}}$ is the set of conjugates of $\pi$ over $\bq$ and $\Phi_{B_{\pi}}\sub \qbar^*$ is the group generated by $\mc{W}_{B_{\pi}}$ . 

Note that a base change of $B_{\pi}$ to the degree $N$ extension  $\bF _{q^N}$ has the effect of replacing $\mc{W}_{B_{\pi}}$ by $\mc{W}_{B_{\pi}}^N := \{\al^N:\al\in \mc{W}_{B_{\pi}}\} = \{\sigma(\pi^N):\sigma\in \absq \}$. Furthermore, $B_{\pi}$ is \textit{absolutely} simple if and only if $\bq(\pi^N) = \bq(\pi)$ for every $N\geq 1$.
\vspace{0.1cm}
\begin{Def} \normalfont An \textit{admissible} Newton polygon of length $2g$ is a sequence $0\leq \lamb_1\leq \cdots \leq \lamb_{2g}\leq 1$ of rationals such that:\\
$(i)$ For each slope $\lamb$, $\lamb$ and $1-\lamb$ occur with the same multiplicity.\\
$(ii)$ If $m_{\lamb}$ is the multiplicity of $\lamb$, then $m_{\lamb}\lamb\in \bz$.\end{Def}

It is well-known that the Newton polygon of an abelian variety over a finite field is admissible. Conversely, every admissible polygon occurs as the Newton polygon of an abelian variety over the algebraic closure of the finite field. Using the results obtained in the preceding section, we now determine the possible polygons at the places of good reduction.

\begin{Prop} Let $A$ be an abelian variety of type $\IV(1,3,\mr{Std})$ over a number field $F$. Let $v$ be a place of good reduction that is prime to the discriminant of $D$. Then $A_v$ has one of the following Newton polygons:\\
$\mr{i}$. $9\times 0, 9\times 1$\\
$\mr{ii}$. $18\times \frac{1}{2}$\\
$\mr{iii}$. $6\times 0, 6\times \frac{1}{2}, 6\times 0$\\
$\mr{iv}$. $9\times \frac{1}{3}, 9\times \frac{2}{3}$.\end{Prop}

\begin{prf} Since $v$ is prime to the discriminant of $D$, Proposition 4.2 implies that $A_v$ is isogenous to $B_v^3$ for some abelian variety $B_v$ of dimension $3$. So it suffices to show that $B_v$ does not have the slope sequence $0,4\times \frac{1}{2}, 1$ as its Newton polygon.

Let $\pi_i,\ov{\pi}_i$ $(i=1,2,3)$ be the eigenvalues of the characteristic polynomial of $B_v$. Set $\al_i := \frac{\pi_i}{\sqrt{q}}$. Fix a $p$-adic valuation $\nu$ such that $\nu(q)=1$. The Newton slopes of $B_v$ are then given by $\nu(\pi_i),1-\nu(\pi_i)$ $(i=1,2,3)$.

Fix a maximal torus $T$ in $G_l^{\der}$. Because of the additive relations between the weights of $T$ on the representation $V_l$, we may order the $\al_i$ so that $\al_1\al_2 = \al_3$. So $\nu(\pi_3) = \nu(\pi_1)+\nu(\pi_2)-\frac{1}{2}$ and one can see that the Newton polygon is not $0,4\times \frac{1}{2}, 1$.\end{prf}

We will see that all of these Newton polygons occur for this class of abelian varieties. The idea is to construct special points on the Shimura variety such that the corresponding CM abelian varieties possess such Newton polygons at some places of good reduction. We then show that there exist generic points with reductions isogenous to those of the special points.

\begin{Lem} Let $B_{\pi}$ be an absolutely simple abelian variety over a finite field $k$ of characteristic $p$ with $[D_{\pi}:\bq(\pi)] = d^2$. If $p$ splits completely in $\bq(\pi)$, then\\ 
$\mr{(i)}$ $B_{\pi}$ has all of its Newton slopes in $\frac{1}{d}\bz$.\\
$\mr{(ii)}$ The least common denominator of the Newton slopes is $d$.\end{Lem}

\begin{prf} Note that since $p$ splits completely in $\bq(\pi)$, it also splits completely in the Galois closure $\bqW$.

$\mr{(i)}$: For any place $v$ lying over $p$, the slope $\frac{w(\pi)}{w(|k|)}$ occurs with multiplicity $d[\bq(\pi)_w:\bq_p]$ and since $[\bq(\pi)_w:\bq_p] =1$ for each $v$, it follows that $d\frac{w(\pi)}{w(|k|)}\in \bz$.

$\mr{(ii)}$: The Hasse invariants of $D$ are given by $[\bq(\pi)_w:\bq_p]\frac{w(\pi)}{w(|k|)} = \frac{w(\pi)}{w(|k|)}$ where $v$ runs through the places of $\bq(\WA)$ lying over $p$. Since $[D_{\pi}]$ is of order $d$ in $\Br(K)$, the least common denominator of its Hasse invariants is $d$.\end{prf}

\begin{Lem} Let $B_{\pi}$ be an absolutely simple abelian variety over a finite field $k$ with\\ $[D_{\pi}:\bq(\pi)] = d^2$. Suppose the extension $\bq(\W_{B_{\pi}})/\bq$ is such that the decomposition group of a prime of $\bq(\W_{B_{\pi}})$ lying over $p$ is a normal subgroup of $\Gal(\bq(\W_{B_{\pi}})/\bq)$. Then:\\ 
$\mr{(i)}$ $B_{\pi}$ has all of its Newton slopes in $\frac{1}{d}\bz$.\\
$\mr{(ii)}$ The least common denominator of the Newton slopes is $d$.\end{Lem}
\vspace{-0.15cm}
\noindent\underline{Remark} In particular, this condition is fulfilled when the center $\bq(\pi)$ of $\End^0(B_{\pi})$ is an abelian extension over $\bq$. 

\begin{prf}  Because of the preceding lemma, it suffices to show that $p$ splits completely in $\bq(\pi)/\bq$. Let $\mf{p}$ be a prime of $\bq(\W_{B_{\pi}})$ lying over $p$. We may assume without loss of generality that $p$ does not split completely in $\bq(\W_{B_{\pi}})$ and hence, the decomposition group of any prime lying over $p$ is non-trivial.
	
Choose a non-trivial element $\tau$ in the decomposition group $D_{\mf{p}}$. Since $\pi$ and its conjugates generate $\bq(\W_{B_{\pi}})$, $\tau$ cannot fix every conjugate of $\pi$. So we may assume without loss of generality that $\tau(\pi)\neq \pi$. Furthermore, the hypothesis implies that $D_{\mf{p}}$ is normal. So $\tau$ lies in the decomposition group of every prime of $\bqW$ lying over $p$ and hence, $\tau$ preserves the ideal $(\pi)$ in $\bq(\W_{B_{\pi}})$. So the algebraic number $\frac{\tau(\pi)}{\pi}$ has the valuation $0$ under every non-archimedean valuation and hence, must be a root of unity. Write $\tau(\pi) = \zeta_{_N}\pi$. Then the characteristic polynomial of $B\times_{\bF _{q}} \bF _{q^N}$ has the zero $\pi^N$ with multiplicity $\geq 2$, a contradiction since $A_v$ is \textit{absolutely} simple. So $p$ splits completely in $\bq(\W_{B_{\pi}})/\bq$. \end{prf}

\begin{Prop} Let $K$ be an imaginary quadratic field, $L^+$ a cubic field Galois over $\bq$ and $L=L^+K$. Let $\{\phi_1,\phi_2,\phi_3\}$ be a CM type for $L$. An abelian variety $X$ with CM type $(L,\{\phi_1,\phi_2,\phi_3\})$ has the following possibilities for the Newton polygon at a place $v$:\\
$1$. If $p:=\charac(v)$ is inert in $L/\bq$, then $X_v$ is supersingular.\\
$2$. If $p$ is inert in $L^+/\bq$, but splits in $K/\bq$, then $X_v$ has slopes $3\times \frac{1}{3},3\times \frac{2}{3}$.\\
$3$. If $p$ splits completely in $L/\bq$, then $X_v$ is ordinary. \end{Prop}

\begin{prf} Basic complex multiplication. \end{prf}


Note that in the case of non-ordinary reduction, lemma 5.3 implies that $\End^0(X_v)$ is a non-commutative division algebra with $L\hra \End^0(X_v)$. If $X_v$ has slopes $3\times \frac{1}{3},3\times \frac{2}{3}$, then $\End^0(X_v)$ is the $9$-dimensional central division algebra over $K$ such that \vspace{-0.15cm}$$\inv_{\mfp}(\End^0(X_v)) = \frac{1}{3},\;\;\; \inv_{\ov{\mfp}}(\End^0(X_v)) = \frac{2}{3}$$ and is unramified everywhere else. In particular, the reduction $X_v$ is a simple three-fold with $p$-rank zero.


\begin{Prop} Let $E$ be a CM field of degree $2g$ Galois over $\bq$. Let $A$ be an abelian variety over a number field $F$ with CM by $E$. For any place $v$ of simple reduction, if $\End^0(X_v)$ is a $d^2$-dimensional central division algebra over a CM field $E$ of degree $2e$, then $d|g$ and the Newton polygon of $A_v$ has slopes in $\frac{1}{d}\bz$.\end{Prop}

\begin{prf} Note that \vspace{-0.15cm}$$ [K:\bq] = 2g = 2e\cdot d\;\; \mr{(Honda-Tate)}$$ and hence, $K$ is a maximal subfield of $D$. The least common denominator of the Hasse invariants of $D$ is $d$. Hence, the Newton slopes of $A_v$ lie in $\frac{1}{d}\bz$.\end{prf}

\begin{center}
\section{\fontsize{11}{11}\selectfont Existence results}
\end{center}
\vspace{-0.1cm}

Let $A$ be an abelian variety of type $\IV(1,d,\wedge^r(\mr{Std}))$. The reductive group $\MT_A$ yields a Shimura variety of PEL type. For the existence of abelian varieties with the reductions described in the preceding section, we use the existence of special points where the abelian variety has such reductions and apply Proposition 6.2.

\begin{Prop} Let $A$ be an abelian variety of type $\IV(1,d,\wedge^r(\mr{Std}))$. The Shimura variety $(\MT_A,X)$ has reflex field $K$.\end{Prop}
\begin{prf} It suffices to show that the intersection of all fields splitting $\MT_A$ is $K$. Now, a number field splits $\MT_A$ if and only if contains a subfield $L$ satisfying the following conditions:\\
-$K\sub L$ \\
-$L$ is a CM field\\
-$L$ splits $D$.

By Proposition 3.4, we can construct two cyclic extensions $L_1^+/\bq$, $L_2^+/\bq$ of degree $d$ such that that $L_1\cap L_2 =\bq$ and $L_1:=KL_1^+$, $L_2:=KL_2^+$ have embeddings in $D$. Since $L_1^+, L_2^+$ are cyclic of odd degree over $\bq$, they are totally real and hence, $L_1,L_2$ are CM fields. So $E(\MT_A,X)\sub L_1\cap L_2 = K$.\end{prf}


\begin{Def}\normalfont For a Shimura variety $(G,X)$ of Hodge type, we say a point $x$ is \textit{generic} if the Mumford-Tate group of $x$ is the same as that of $X$ and \textit{special} if the Mumford-Tate group is a torus.\end{Def}

Let $A$ be an abelian variety of type $\IV(1,3,\mr{Std})$ with endomorphism algebra $D$. For brevity, we write $G = \MT_A = D^{\Nm=1}$. Let $\Gam_{\infty}\sub G(\br)$ be the centralizer of $h:\mb{S}\lra G_{\br}$. Then for any sufficiently small torsion-free compact subgroup $\Gam\sub G^{\der}(\bq)$, we may form the connected Shimura variety $M_{\bc}(\bc) = \Gam\bs G(\br)^0/\Gam_{\infty}^0$. This is the set of $\bc$-valued points of a quasi-projective $\bc$-scheme and is a connected component of a moduli space for a class of abelian varieties over $\bc$. Replacing $\Gam$ by a smaller subgroup if necessary, we may assume that $M_{\bc}$ is smooth and there exists a universal polarized abelian scheme $\mc{X}_{\bc}\lra M_{\bc}$.

Let $x$ be a special point on this Shimura variety, meaning that the associated map $h:\mb{S}\lra G_{\br}$ factors through a maximal torus $T_x\sub G$. By [CO12], there exists a CM field $L\sub D$ with $[L:K] =3$ such that $T_x$ is the kernel of the map \vspace{-0.2cm }$$(\Nm_{L/L^+},\Nm_{L/K}): \Res_{L/\bq} \bg_m \lra \Res_{L^+/\bq}(\bg_m)\times\Res_{K/\bq}(\bg_{m,\bq}).$$ Let $L^+$ be the maximal real subfield of $L$ and let $\wti{L}, \wti{L^+}$ be the Galois closures of $L,L^+$ over $\bq$. Then $L=L^+K$, $\wti{L} = \wti{L^+}K$

By the last proposition, $M_{\bc}$ has a quasi-canonical model $M$ over $E(G,X)=K$. Enlarging $E(G,X)$ if necessary, we may assume the abelian scheme $\mc{X}_{\bc}$ descends to an abelian scheme $\mc{X}/M$.

\begin{Prop} Let $A$ be an abelian variety of type $\IV(1,3,\mr{Std})$. Let $(G,X,h)$ be a Shimura variety with generic Mumford-Tate group $G = \MT_A$. Let $\mfp$ be a prime of a field $F$ contaning the reflex field $K$. Let $x\in X(F)$ be a special point corresponding to an abelian variety $A(x)$ with CM. Then there exists a finite extension $F'/F$, a prime $\mfp'$ lying over $\mfp$ and a generic point $y\in X(F')$ corresponding to an abelian variety $A(y)$ such that $A(x)$ and $A(y)$ have isogenous reductions at $\mfp'$.\end{Prop}

\begin{prf} Enlarging the field if necessary, we may assume $X$ has good reduction at $v$. Let $\mc{O}_v$ be the localization of $\mcO_F$ at $v$. We may replace the compact subgroup $C\sub G(\mb{A}_f)$ by a compact subgroup so that it defines a level structure that is prime to $p:=\charac(v)$ and such that there exists a finite map $M\lra \mc{A}_9$ where $\mc{A}_9$ is a fine moduli scheme of $9$-dimensional abelian schemes with level structure prime to $p$. We replace $M$ by the absolutely irreducible component containing $x_F$. 

Let $\mc{M}/\mc{O}_v$ be the Zariski closure of $M$ in $\mc{A}_9$. The image in $\mc{A}_{9}(F)$ of $x_F\in M(F)$ extends to $x\in \mc{M}(\mc{O}_{v})$. We replace $\mc{M}$ by an affine open	subset containing $x$ and choose $\pi:\mc{M}\lra \mb{A}_{\mc{O}_F}^1$ to be a map such that $\pi_{k_v}$ is non-constant on every irreducible component of $\mc{M}_{k_v}$ and $\pi_F:\mc{M}_F\lra \mb{A}_{F}^1$ is etale. It follows from ([No95], Proposition 1.7) that there is a thin subset $S\sub\mb{A}^1(F)$ such that each $y_F\in \mb{A}^1(F) \bs S$ corresponds to an abelian variety with Mumford-Tate group $G$. Since $F$ is Hilbertian, it follows that there exist infinitely many points $y'\in \mb{A}^1(\mc{O}_{F})$ with $y'_{k_v} = \pi(x)_{k_v}$ and $y'_F\notin S$.\end{prf}

\begin{center}
\section{\fontsize{11}{11}\selectfont Appendix}
\end{center}
\vspace{-0.2cm}

We construct the abelian varieties described in Definition 1.1. This construction is a slight modification of that in ([Orr15], Section 5.2). Recall that we have an equivalence between the categories of abelian varieties over $\bc$ and polarizable $\bz$-Hodge structures of type $\{(-1,0),(0,-1)\}$ given by $X\mapsto H^1(X(\bc),\bz)$.

Let $d$ be an odd integer and let $r$ be an integer such that $1\leq r\leq d$ and $\gcd(r,d)=1$. Let $K$ be a an imaginary quadratic field and let $\wti{D}$ be a $d^2$-dimensional central division algebra over $K$ with an involution $*$ of the second kind. We obtain the $\bq$-algebraic groups \vspace{-0.1cm}$$ H(R):= \{\al\in (\wti{D}\otimes_{\bq}R)^{\times}: \al\al^*=1\},\;\;\; G(R):= \{\al\in (\wti{D}\otimes_{\bq}R)^{\times}: \al\al^*\in R\}$$ which are $\bq$-forms of $\SL_d$ and $\GL_d$ respectively.

We view $\wti{D}$ as a $K$-irreducible representation of $H_K$ with highest weight $\ovm_1$. Let $D$ be the unique central division algebra over $K$ such that $[D]^{-1} = [\wti{D}]^r$ in $\Br(K)$. Since $r$ and $d$ are relatively prime, $[D]$ is of order $d$ and hence, $[D:K] = d^2$. 

Let $\wti{\rho}$ be the $K$-irreducible representation of $H_K$ with highest weight $\ovm_r$. Then $\wti{\rho}$ has endomorphism algebra $D$ and $\wti{\rho}\otimes_K \qbar$ is the sum of $d$ copies of the $r$-th exterior product of the standard representation of $\SL_d$. So $\dim_K \wti{\rho} = d {d\choose r}$.

If $\lamb I$ is a scalar matrix in $H(\bc)$, then $\wti{\rho}_{_{\bc}}(\lamb I)$ is multiplication by $\lamb^r$ and we may extend $\wti{\rho}$ to a representation of $G_K$ by letting each scalar matrix act as multiplication by $\lamb^r$.

Set $\rho = \Res_{K/\bq} \wti{\rho}$. This is a $\bq$-irreducible representation of $G$ of dimension $2d {d\choose r}$. We have $\ker(\rho) = \mu_r$. So the resulting representation of $G/\mu_r$ is faithful.

Now, $H_{\br}$ is the unitary group of the Hermitian form of signature $(1,d-1)$. We define $\phi':\bc^*\lra G(\br)$ as follows. Let $\phi'(z)$ act as $\frac{z^r}{\ov{z}^{r-1}}$ on the space where $h$ is positive definite and as $\ov{z}$ on the space where $h$ is negative definite. Then $\rho\circ	\phi'$ has weights $z^r$ and $\ov{z}^r$ and since $\rho$ is a faithful representation of $G/\mu_r$, there is a homomorphism $\phi:\bc^*\lra (G/\mu_r)(\br)$ whose $r$-th power is $\phi'$. So $(G/\mu_r,\rho,\phi)$ defines a $\bq$-Hodge structure of type $\{(-1,0),(0,-1)\}$ and the Hermitian form induces a polarization of the Hodge structure.

An abelian variety $X$ over $\bc$ in the corresponding isogeny class has Mumford-Tate group $G/\mu_r$, endomorphism algebra $D$ and dimension $d {d\choose r}$. We may then choose a $\qbar$-valued generic point on the Shimura variety attached to $\MT_X$. This will correspond to an abelian variety over $\qbar$ that fulfills the conditions of Definition 1.1.

\vspace{0.4cm}

\noindent\textbf{Acknowledgments:} The author is grateful to Martin Orr for patiently answering his questions.

\vspace{0.2cm}

\begin{center}
	\textbf{References} 
\end{center}
\vspace{-0.2cm}
\footnotesize 

\noindent [Bog80] Fedor Alekseivich Bogomolov, Sur l'algebricite des representations l-adiques, C. R. Acad. Sci. Paris Ser.
A-B 290 (1980), no. 15, A701A703. MR574307

\noindent[Chi92] W. Chi, \textit{On the $l$-adic representations attached to simple abelian varieties of type $\IV$}, Bull. Aust. Math Soc. 1992

\noindent[CO12] C.L.Chai and F. Oort, \textit{Abelian varieties isogenous to a Jacobian}, Ann. of Math. (2) 176(1) (2012), 58635

\noindent[CCO13] B. Conrad, C. Chai, F. Oort, \textit{Complex multiplication and lifting problems}, AMS series Mathematical Surveys and Monographs, Vol. 195. AMS 2013

\noindent[Fal86] G. Faltings, \textit{Finiteness theorems for abelian varieties over number fields}, Arithmetic Geometry, 1986

\noindent [Gro66] A. Grothendieck, Elements de geometrie algebrique. IV. Etude locale des schemas et des morphismes de ´
schemas. III, Inst. Hautes Etudes Sci. Publ. Math. 28 (1966), 255. MR0217086 

\noindent[Jac80] N. Jacobson, \textit{Basic Algebra Vol. II}, 1980

\noindent[KMRT98] Knus, Merkujev, Rost, Tignol, \textit{The book of involutions}, American Mathematical Society Col. Publications, 44, AMS 1998

\noindent[Per12] K.Madapusi Pera, \textit{Toroidal compactifications of integral models of Shimura varieties of Hodge type}, preprint 2015

\noindent[LP97] M. Larsen and R. Pink , \textit{A connectedness criterion for l-adic Galois representations}, Israel J. Math. 97 (1997).

\noindent [Mum74] D. Mumford, \textit{Abelian Varieties}, 2nd ed., Oxford University Press, 1974

\noindent [MFK94] D. Mumford, J. Fogarty, and F. Kirwan, \textit{Geometric invariant theory}

\noindent[MP08] V.K. Murty and V. Patankar, \textit{Splitting of abelian varieties}, Int. Math. Res. Not. IMRN 12, 2008

\noindent[No95] R. Noot, \textit{Abelian varieties- Galois representations and properties of ordinary reduction}, Comp. Math. 97

\noindent[No09] \underline{\;\;\;\;\;\;\;} , \textit{Classe de conjugaison du Frobenius d'une variete abelienne sur un corps de nombres}, J. Lond. Math. Soc. (2) 79 (2009), no. 1.

\noindent[Orr15] M. Orr, \textit{Lower bounds for ranks of Mumford-Tate groups}, Bull. Soc. Math. France, 143(2), 2015

\noindent[Pin98] R. Pink, \textit{l-adic algebraic monodromy groups, co-characters and the Mumford-Tate conjecture}, J. Reine Angew. Math. 495 (1998)

\noindent[Pin04] \underline{\;\;\;\;\;\;\;},\textit{On the order of the reduction of a point on an abelian variety}, Math. Ann. 330 (2004), no. 2, 275–291

\noindent[Ser79] J.P. Serre, \textit{Groupes alg´ebriques associ´es aux modules de Hodge-Tate, Journ´ees de G´eom´etrie Alg´ebrique de Rennes} (Rennes, 1978), Vol. III, 1979, pp. 155–188.

\noindent[ST68] J.P. Serre and J.Tate, \textit{Good reduction of abelian varieties}, Ann. of Math. (2), 88, 1968

\noindent[Th16] S. Thakur, \textit{Reductions of abelian varieties of generalized Mumford type}, Preprint (Arxiv)

\noindent[Zyw14] D. Zywina, \textit{The splitting of reductions of an abelian variety}, IMRN, vol. 2014, No. 18.

\vspace{0.4cm}

\noindent Email: stevethakur01@gmail.com

\end{document}